\documentclass{article}
\usepackage{amssymb,amsmath}
\date{}
\author{}
\title{}
\def\dps{\displaystyle}
\def\qed{\hspace{\fill}$\Box$}
\title{Left invariant semi Riemannian metrics on quadratic Lie groups}
\author{S. Bromberg  \and A. Medina}
\newtheorem{theorem}{Theorem}
\newtheorem{pr}{Proposition}
\newtheorem{lm}{Lemma}
\newtheorem{df}{Definition}
\newtheorem{cor}{Corollary}
\newtheorem{rmk}{Remark}
\date{}

\begin{document}

\noindent {\centerline{{\Large{Left invariant semi Riemannian
metrics on }}}}

{\centerline{{\Large{quadratic Lie groups}}}}

\vspace{1in} \noindent Shirley Bromberg\\Departamento de Matem\'aticas\\
Universidad Aut\'onoma Metropolitana-Iztapalapa\\ M\'exico, D.F. M\'exico\\
stbs@xanum.uam.mx

\vspace{.3in} \noindent Alberto Medina\\D\'epartement des Math\'ematiques, C.C. 051\\
Universit\'e Montpellier 2 UMR CNRS 5149\\ Place E. Bataillon, 34095
Montpellier cedex 5, France \\ medina@math.univ-montp2.fr

\vspace{.3in}
\begin{abstract}
To determine the Lie groups that admit a flat (eventually complete)
left invariant semi-Riemannian metric is an open and difficult
problem. The main aim of this paper is the study of the flatness of
left invariant semi Riemannian metrics on quadratic Lie groups i.e.
Lie groups endowed with a bi-invariant semi Riemannian  metric. We
give a useful necessary and sufficient condition that guaranties the
flatness of a left invariant semi Riemannian metric defined on a
quadratic Lie group. All these semi Riemannian metrics are complete.
We show that there are no Riemannian or Lorentzian flat left
invariant metrics on non Abelian quadratic Lie groups, and that
every quadratic 3 step nilpotent Lie group admits a flat left
invariant semi Riemannian metric. The case of quadratic 2 step
nilpotent Lie groups is also addressed.

\end{abstract}

\vspace{.3in}\noindent {\bf{Key words:}}

Left invariant semi-Riemannian metrics, flat semi Riemannian
metrics, geodesically complete manifolds, quadratic Lie groups,
Jacobi fields.

\section*{Introduction}

 This article outlines some facts known by the authors about the semi
 Riemannian geometry of a Lie group provided with a semi
 Riemannian metric invariant under left translations.

\newpage
Contents

1. General results about left invariant semi Riemannian metrics on Lie groups

2. Quadratic Lie groups

3. Jacobi fields on quadratic Lie groups, conjugate points.

4. On flat left invariant semi Riemannian metrics on quadratic Lie
groups.

5. The nilpotent quadratic case.

When the relations between curvature or (geodesically) completeness
of a semi Riemannian metric and another topological or geometrical
properties are studied it is very useful to have many examples. This
paper describes a rich collection of examples which are obtained by
providing a Lie group $G$  with a semi Riemannian metric invariant
under left translations. It is well known that every left Riemannian
metric is complete and in \cite{kn:Mi} Milnor described the Lie
groups with flat left invariant Riemannian metrics. By contrast, the
study of completeness and/or flatness of a non definite metric is in
general very difficult. Even in the 3 dimensional non unimodular
case, there is not in the literature  a necessary and sufficient
condition that guaranties  the completeness of a left invariant
Lorentzian metric. When the 3 dimensional Lie group is unimodular,
the completeness of a left invariant Lorentzian metric is equivalent
to the completeness of the geodesics of light type
(\cite{kn:bmsigma}).

Our class of examples can be enlarged substantially, with no extra
effort, as follows. If $\Gamma $ is any discrete subgroup of $G,$ then
a left invariant semi Riemannian metric on $G$  gives rise to a
metric on the quotient space  $\Gamma \backslash G $  with identical
properties of curvature and (in)completeness. The case where $\Gamma
\backslash G $ is compact is of particular interest.

The first section will survey general old and new results on left
invariant semi Riemannian metrics on Lie groups.

The principal and new result (Theorem \ref{theorem:ACH}) gives
necessary and sufficient conditions for the flatness of a left
invariant semi Riemannian metric on unimodular Lie groups. Under
these conditions flatness implies completeness.

In section 2,  necessary and sufficient conditions that guarantees
the existence  of  bi-invariant semi Riemannian metrics on Lie
groups are given. These groups, called quadratic or orthogonal Lie
groups, are the central objets of our study. Section 3 is devoted to
the Jacobi vector fields corresponding to a left invariant semi
Riemannian metric on  Lie groups and on quadratic Lie groups in
particular. The equation that defines the reflection on the Lie
algebra of $G$ of a such vector  field is particularly simple  when
the metric is bi-invariant. The reflections of the Jacobi vector
fields corresponding to the Lorentzian bi-invariant metrics on the
oscillator Lie groups are determined.

Theorem \ref{theorem:ACHQ}, Theorem \ref{theorem:ind>2}, Theorem
\ref{theorem:dim4}, and Theorem \ref{theorem:fder3} are the main
results of section 4. The first one specializes Theorem
\ref{theorem:ACH} to the case of quadratic Lie groups. One of the
consequences of Theorem \ref{theorem:ind>2} is that every left
invariant Lorentzian metric on a non Abelian quadratic Lie group is
non  flat. Theorem \ref{theorem:dim4} shows the non existence of
flat left invariant semi Riemannian metric on any indecomposable
quadratic Lie group of dimension 4. The case of left invariant semi
Riemannian metrics on quadratic nilpotent Lie groups is also
treated. Theorem \ref{theorem:fder3} shows that every 3 step
nilpotent Lie group admits a flat left invariant connection given by
an invertible $f$-derivation. This connection is the Levi Civita
connection of a semi Riemannian metric if the group is quadratic
(Theorem \ref{theorem:fder}). A left semi Riemannian metric on a
nilpotent quadratic Lie group $ (G, k ) $ defined by a $k$ symmetric
linear isomorphism $u$ is complete when $u$ preserves the descending
central sequence of the Lie algebra ${\cal{G}}$ ( Proposition 9 ).
If $( G ,k )$ is 2 step nilpotent and its corank is $ 0 $ then  $G$
admits many non isometric flat left invariant semi Riemannian
metrics. If ${\mathrm{dim}}\, G > 8$  there are infinitely many non
isometric such metrics (Theorem \ref{theorem:2n}).

The following result is an important and final remark concerning the classical
or generalized solutions  of the Yang-Baxter equation on quadratic
Lie groups and the relations with  left invariant semi Riemannian
metrics.

\begin{theorem}[\cite{kn:Bou-M}]
Every solution of the classical Yang-Baxter equation on a quadratic
Lie group induces a flat left invariant semi Riemannian metric on
the dual Lie groups associated to the solution. Furthermore a
solution of the generalized Yang-Baxter equation determines a left
invariant semi Riemannian metric such that the covariant derivative
of its curvature tensor vanishes.
\end{theorem}

\section{General results about left invariant semi
Riemannian metrics on Lie groups}

Let $G$ be a Lie group, $\varepsilon$ the unit element in $G.$ A non
degenerate symmetric bilinear form $\langle,\rangle$ on
${\cal{G}}:=G_\varepsilon$ defines a left invariant semi Riemannian
metric on $G$ given by the formula
$$
\langle v_\sigma, w_\sigma\rangle_\sigma:=
\langle(\mathrm{L}_{\sigma^{-1}})_{*,\sigma} v_\sigma,
({\mathrm{L}}_{\sigma^{-1}})_{*,\sigma}w_\sigma\rangle,\quad
\sigma\in G, v_\sigma,w_\sigma\in G_\sigma
$$
where ${\mathrm{L}}_\sigma: \tau\mapsto \sigma\tau,$ and conversely.

The Levi-Civita connection $\nabla$ of a semi Riemannian left
invariant metric is left invariant, and defines a product on the Lie
algebra given by the formula
$$
xy^+:=\nabla_{x^+}y^+,
$$
where $x^+$ stands for the left invariant vector field with
infinitesimal generator $x\in{\cal{G}}.$ This product, called the
{\bf{Levi-Civita product,}} verifies the Koszul formula
$$
2\langle xy,z\rangle = \langle [x,y],z\rangle -\langle
[y,z],x\rangle+ \langle [z,x],y \rangle.
$$
By means of
$$
x(t):=\left({\mathrm{L}}_{\sigma(t)^{-1}}\right)_{*,\sigma(t)}\sigma'(t),
$$
the equation for the geodesics of the semi Riemannian metric
becomes, in the Lie algebra,

\begin{equation}\label{eq:geo}
\dot{x}=-xx.
\end{equation}

Since the Levi-Civita connection is torsion free, the Levi-Civita
product satisfies
$$xy-yx=[x,y].$$

Moreover the Koszul formula implies that the map $L_x: y\mapsto xy$
is $\langle\, ,\, \rangle$ skew symmetric.

Many features of the geometry of left invariant semi Riemannian
metrics on Lie groups can be studied in the Lie algebra.

A semi Riemannian metric is called {\bf{complete}} when its
geodesics are defined for all $t\in\mathbb{R}.$ Notice that a left
invariant semi Riemannian metric is complete if and only if the
solutions of equation (\ref{eq:geo})  are defined for all values of
the parameter.

The exponential map associated to a semi Riemannian metric with base
point $\sigma\in G$ is denoted by ${\mathrm{Exp}}_\sigma.$ This map
is defined on all of $G_\sigma$ for all $\sigma$ whenever the semi
Riemannian metric is complete. Notice that in general
${\mathrm{Exp}}_\varepsilon$ differs from the exponential map in Lie
theory (see remark \ref{rmk:2} bellow).

\begin{df}
A semi Riemannian metric is called {\bf{flat}} if the curvature
tensor vanishes.
\end{df}

The Levi-Civita product for a flat left invariant semi Riemannian
metric on a Lie Group $G$ is a left symmetric product on
${\cal{G}},$ compatible with the Lie bracket, that is
$$(xy)z-x(yz)=(yx)z-y(xz),$$
and $$xy-yx=[x,y].$$

As a partial converse we have that a left symmetric product
compatible with the Lie bracket induces a flat left invariant
connection on $G.$

\vspace{.1in}

The existence of a flat left invariant  metric on a Lie group
imposes serious restrictions on the group as the following result
shows

\begin{theorem}[Theorem 1.5 \cite{kn:Mi}]
A Lie group has a left invariant flat Riemannian metric if and only if its
Lie algebra decomposes as a semidirect product of an Abelian Lie
algebra with an Abelian ideal, the Abelian algebra acting on the
Abelian ideal by infinitesimal isometries.
\end{theorem}

\vspace{.1in}

The existence of a flat left invariant  metric on a quadratic Lie
group imposes even more restrictive conditions on the group. In fact
under this hypothesis the group is Abelian (see proposition
\ref{pr:fq}). In the same line of thought Proposition 7 states that
on non Abelian quadratic Lie groups there are no flat left invariant
Lorentzian metrics.

\vspace{.1in} In what follows, the following result will be useful
\begin{theorem}\label{theorem:uni}
A flat left invariant semi Riemannian metric is complete if and only
if the Lie group is unimodular.
\end{theorem}

For the proof see \cite{kn:au-m}.

\vspace{.1in} In \cite{kn:g-m} the simply connected Lie groups with
flat, complete left invariant Lorentz metrics are characterized. The
nilpotent case was treated alternaternatively by means of the double
extension in \cite{kn:au-m}.

The Jacobi fields measure the variation of geodesics: If $t\mapsto
\tau (t)$ is a geodesic, the vector field $t\mapsto Y(t)$ defined on
the curve $\tau$ is a {\bf{Jacobi vector field}} provided that it
satisfies the second order differential equation
\begin{equation}\label{eq:Jacobi}
\frac{D^2Y}{dt^2}= {\mathrm{R}}_{Y\tau'}(\tau')
\end{equation}
where $DY/dt$ stands for the affine covariant derivative of
$(G,\nabla)$ of the vector field $Y$ on the curve $\tau$ and
${\mathrm{R}}$ is the curvature tensor.

\vspace{.2in}

Hence if the semi Riemannian metric is flat, then a vector field $Y$
on a geodesic is a Jacobi field if and only if the vector field
$\dfrac{DY}{dt}$ is parallel along the geodesic, that is if and only
if $\dfrac{D^2Y}{dt^2}=0.$

Then a necessary condition for a left invariant semi Riemannian
metric to be flat is that the second covariant derivative of any
right invariant vector field along every geodesic vanishes.

\vspace{.1in} Notice that every right invariant vector field is a
Jacobi vector field along any geodesic, because it is a Killing
vector field (~\cite{kn:K-N}).

\begin{df} Let
$\tau:[a,b]\to G$ a geodesic. The points $\tau(a),$ $\tau(b)$ are
called {\bf{conjugate points}} if there is a Jacobi field $Y$ on
$\tau$ such that $Y(a)=Y(b)=0.$
\end{df}

\begin{pr}
Let $\tau:[a,b]\to G$ be a geodesic. Then $\tau(a)$ and $\tau(b)$
are conjugate if and only if the rank of ${\mathrm{Exp}}_{\tau(a)}$
at $(b-a)\tau'(a)$ is less than dim $G$.
\end{pr}

\begin{lm}\label{lm:isom}
Let $(M,\langle\,,\,\rangle)$ a flat semi Riemannian  manifold. Let
$U$ be a connected neighborhood of $0\in M_\sigma$ where
${\mathrm{Exp}}_\sigma$ is defined. If the semi Riemannian metric is
flat then ${\mathrm{Exp}}_\sigma$ is a local isometry.
\end{lm}

\vspace{.1in}\noindent{\bf{Proof.}} Let $x\in M_\sigma,$
$v,w\in(M_\sigma)_x \approx M_\sigma.$ Let $J_v,J_w$ the unique
Jacobi fields along a geodesic $\tau:[0,1]\to G,$ $\tau(0)=\sigma,
\tau(1)=\rho$ such that
$$J_v(0)=J_w(0)=\sigma,
\dfrac{DJ_v}{dt}(0)=v,\dfrac{DJ_w}{dt}(0)=w.$$
The derivatives of the map
$$
\varphi(t):=\langle J_v(t),J_w(t)\rangle.
$$
are
$$
 \varphi'(t)=\langle \dfrac{DJ_v}{dt}(t),J_w(t)\rangle +  \langle
 J_v(t),\dfrac{DJ_w}{dt}(t)\rangle,
$$
$$
\varphi''(t)=2\langle\dfrac{DJ_v}{dt},\dfrac{DJ_w}{dt}\rangle, \varphi'''(t)= 0
$$
since $J_v,J_w$ are parallel along $\tau.$ Hence
$\varphi''(t)=2\langle v, w\rangle,$ $\varphi'(t)=2t\,\langle v,
w\rangle,$ $\varphi(t)=t^2\,\langle v, w\rangle,$ and
$$
\langle d{\mathrm{Exp}}_\sigma(0)v,
d{\mathrm{Exp}}_\sigma(0)w\rangle = \varphi (1)= \langle v,
w\rangle.
$$

\qed

As a corollary we have that a complete semi Riemannian flat metric
has no conjugate points. Furthermore

\begin{lm}
Let $(M, \langle\,,\,\rangle)$ be a flat semi Riemannian manifold.
If ${\mathrm{Exp}}_p$ is defined for all $v\in M_p,$ then
${\mathrm{Exp}}_p:(M_p,\langle\!\langle\,,\,\rangle\!\rangle) \to
(M, \langle\,,\,\rangle)$ is a semi Riemannian covering, where
$\langle\!\langle\,,\,\rangle\!\rangle$ is the affine metric induced
by $\langle\,,\,\rangle_p.$
\end{lm}

\vspace{.1in}\noindent {\bf{Proof.}} We have to show that
${\mathrm{Exp}}_p$ has the lifting property for geodesics. Let
$\tau:[0,1]\to M$ a geodesic and $x_0\in M_p$ such that
${\mathrm{Exp}}_\sigma(x_0)=\tau (0).$ By lemma \ref{lm:isom}, there
are neighborhoods $U$ and $V$ of $x_0$ in $M_p$ and $\tau(0)$ in $M$
such that ${\mathrm{Exp}}_p$ defined on $U$ onto $V$ is an isometry.
If $t$ satisfies $\tau ([0,t))\subset V,$ then $ s\mapsto
{\mathrm{Exp}}_p^{-1}\tau (s)$ is a geodesic in $M_p.$ By hypothesis
this geodesic is defined in $[0,1]$ and it is a lifting of $\sigma.$
The conclusion follows from Theorem 28.7 in \cite{kn:ON}.

\qed

\vspace{.1in}
The following theorem puts together some scattered results

\begin{theorem}\label{theorem:ACH}
Let $G$ be a connected unimodular Lie group and
$\langle\,,\,\rangle$ a left invariant semi Riemannian metric on
$G.$ Then the following assertions are equivalent:

\vspace{.1in}\noindent i)  $\langle\,,\,\rangle$ is flat and complete.

\noindent ii) ${\mathrm{Exp}}_\varepsilon$  is a local isometry
(hence, for every $\sigma$ in $G,$ ${\mathrm{Exp}}_\sigma$ is a
local isometry).

\noindent iii) $\langle\,,\,\rangle$ is flat.

\vspace{.1in}\noindent In any case $G$ is solvable, and
$(G,\langle\,,\,\rangle)$ has no conjugate points.
\end{theorem}

\vspace{.1in}\noindent{\bf{Proof.}} By Theorem \ref{theorem:uni} a
flat left invariant semi Riemannian metric defined on an unimodular
Lie group is complete. The hypothesis of $\langle\,,\,\rangle$ being
flat implies that $G$ is locally symmetric and that
${\mathrm{Exp}}_q$ is a local isometry that has the lifting property
for geodesics. Hence it is a semi Riemannian covering.

To prove that $G$ is solvable notice that  the hypothesis imply that
the Levi-Civita connection defined by the semi Riemannian metric is
a left invariant affine structure. Then there is a representation
$\theta$ of $G$ by affine transformations of ${\cal{G}}$ with an
open orbit and discrete isotropy (see \cite{kn:M-JDG}). The action
of $G$ on ${\cal{G}}$ induced by $\theta$ is transitive because the
metric is complete. Hence the restriction of the representation to a
Levi subalgebra is completely reducible. This contradiction implies
the solvability. \qed

\vspace{.15in}\noindent{\bf{Example.}} Let
$G=\mathbb{R}^2\rtimes{\mathrm{SO}}(\mathbb{R}^2)$ be the connected
component of the unit element of the group of rigid motions of the
plane. The product on $G$ is given by
$$
(x,y,\alpha)(x',y', \beta)= (x+x'\cos\alpha-y'\sin\alpha,
y+x'\sin\alpha+y'\cos\alpha, \alpha +\beta).
$$
Let ${\cal{G}}= {\mathrm{Span}}\{e_1,e_2\}\rtimes \mathbb{R}\,e_3, $
where ${\mathrm{e_1,e_2}}$ is an Abelian Lie ideal  and
$[e_3,e_1]=e_2,\; [e_3,e_2]=-e_1.$ Define a left invariant semi
Riemannian metric by the Lorentzian quadratic form on ${\cal{G}}:$
$$
q( x_1e_1+x_2e_2+x_3e_3)=x_1^2+x_2^2-x_3^2.
$$
Some straightforward calculations show that
$$
L_{e_1}=L_{e_2}= 0 \quad\quad L_{e_3}={\mathrm{ad}}_{e_3}.
$$
Then $L_{[x,y]}=0$ and $L_{x}L_{y}=L_{y}L_{x}.$ Hence the Lorentzian
metric is flat. Equation (\ref{eq:geo})  is in this case
\begin{eqnarray*}
\dot{x_1}&=& x_2x_3\\
\dot{x_2}&=& -x_1x_3\\
\dot{x_3}&=& 0
\end{eqnarray*}
The solution to this equation with initial condition $(x_1,x_2,x_3)$
is
$$
x(t)=
\left(x_1\cos(x_3t)-x_2\sin(x_3t),x_1\sin(x_3t)+x_2\cos(x_3t),x_3\right).
$$
The geodesic on $G$ starting at $\varepsilon$ with initial velocity
$(x_1,x_2,x_3)$ is for $x_3=0$
$$
\gamma(t)= (tx_1,tx_2,0),
$$
and when $x_3\ne 0:$
$$
\gamma(t)=
(-\frac{x_2}{2x_3}+\frac{x_2}{2x_3}\cos(2x_3t)+\frac{x_1}{2x_3}\sin(2x_3t),
\frac{x_1}{2x_3}-\frac{x_1}{2x_3}\cos(2x_3t)+\frac{x_2}{2x_3}\sin(2x_3t),x_3t).
$$
Hence the exponential map based at $\varepsilon$ is, for $x_3=0,$
$$
{\mathrm{Exp}}_\varepsilon(x_1,x_2,x_3)= (x_1,x_2,0)
$$
and for $x_3\ne 0,$
$$(-\frac{x_2}{2x_3}+\frac{x_2}{2x_3}\cos(2x_3)+\frac{x_1}{2x_3}\sin(2x_3),
\frac{x_1}{2x_3}-\frac{x_1}{2x_3}\cos(2x_3)+\frac{x_2}{2x_3}\sin(2x_3),x_3). $$
Hence $
{\mathrm{Exp}}_\varepsilon$  is a global isometry.

\section{Quadratic Lie groups}

\begin{df}
A Lie group $G$ with a bi-invariant semi Riemannian metric $k$ is
called {\bf{orthogonal}} or {\bf{quadratic}} Lie group. The pair
$({\cal{G}}, k),$ where ${\cal{G}}$ is the corresponding Lie algebra
and $k$ is the restriction of $k$ to ${\cal{G}},$ is called
{\bf{orthogonal}} or {\bf{quadratic}} Lie algebra.
\end{df}

Let $({\cal{G}}, k)$ be a quadratic Lie algebra. Then $k$ is a non
degenerate quadratic form and ${\mathrm{ad}}_x$ is $k$ skew
symmetric for all $x\in{\cal{G}}:$
$$\label{eq:quad}
k(\mathrm{ad}_x\,y,z)+k(y,\mathrm{ad}_x\,z)=0.
$$

For every left invariant semi Riemannian metric $\langle\, ,\rangle$
on $G$ there is a  $k$ symmetric isomorphism $u$ of the vector space
underlying ${\cal{G}}$  such that, for all $x,y \in{\cal{G}}$
$$
\langle x , y\rangle = k(u(x),y).
$$
Equation (\ref{eq:geo})  becomes in this case
$$
u(\dot x)=[u(x),x].
$$

\vspace{.1in} The following propositions characterizes quadratic Lie
groups.

\begin{pr}[\cite{kn:MRCA}]
A Lie group is quadratic if and only if the adjoint and co-adjoint
actions are isomorphic.
\end{pr}

\begin{pr}
The Lie group $G$ is quadratic if and only if the linear Poisson
structure on ${\cal{G}}^*$ given by the Lie bracket of ${\cal{G}}$
has a quadratic, non degenerate Casimir.
\end{pr}

\noindent {\bf{Proof.}} Let  $({\cal{G}},k)$ be an orthonormal Lie
algebra. Denote by $\Phi : {\cal{G}}\to{\cal{G}}^*$ the symmetric
isomorphism $\Phi(x):= k(x, \cdot).$ For $x\in{\cal{G}},$ define
$\hat{x}\in({\cal{G}}^*)^*$ by $\hat{x}(\xi)=\xi(x),$
$\xi\in{\cal{G}}^*.$

Let $f:{\cal{G}}^*\to\mathbb{R}$ be given by
$f(\xi)=\xi(\Phi^{-1}(\xi)).$ Clearly $f$ is a non degenerate
quadratic form, hence its differential is $(d\, f)_\xi(\xi')=
2\xi'\left(\Phi^{-1}(\xi)\right)$ or, equivalently, $(d\,
f)_\xi=2\Phi^{-1}(\xi)^{\wedge}.$ If $x_0=\Phi^{-1}(\xi)$, we get,
using the $\mathrm{ad}$ invariance of $k,$ that
$$
\{f,\hat{x}\}_\xi
:=\xi[2\Phi^{-1}(\xi),x]=\Phi(x_0)[2x_0,x]=k(x_0,[2x_0,x])=0,
$$
for all $x\in{\cal{G}}.$ Therefore $f$ is a  Casimir for the
Lie-Poisson bracket  $\{\cdot,\cdot\}.$

Conversely, let  $f(\alpha)={\bf{b}}(\alpha,\alpha)$ be a Casimir,
${\bf{b}}$ being a quadratic, symmetric and non degenerate quadratic
form. Define $k:{\cal{G}}\times{\cal{G}}\to\mathbb{R}$ by
$k(x,y):=\Psi^{-1}(\hat{x})y,$ where
$\Psi(\alpha):={\bf{b}}(\alpha,\cdot).$ Since $\Psi$ is a symmetric
isomorphism, so is $k .$ Moreover, $f$ being a Casimir, we get that
for all $x\in{\cal{G}},$ $\{f,\hat{x}\}=0,$ that is,
$$
\forall\,\alpha\in{\cal{G}}^*,\quad\forall\, x\in{\cal{G}}, \;
0=\{(df)_\alpha,\hat{x}\}_\alpha=\{\Psi(\alpha), \hat{x}\}_\alpha.
$$
Hence for all $y,z\in{\cal{G}},$
$$k(z,[z,y])=0.$$  Replacing $z$ by $a+b,$  $a,b\in{\cal{G}},$ we
get that
$$
\forall\, a,b,y\in{\cal{G}}, \;\;k(a,[b,y])+k(b,[a,y])=0.
$$
Hence $k$ is a quadratic structure on ${\cal{G}}.$ \qed

\vspace{.2in}

\begin{rmk}\label{rmk:2}
The Levi-Civita product and the curvature of a bi-invariant semi
Riemann\-ian metric on a Lie group are given (resp.) by:
$$
xy=\frac{1}{2}[x,y],\quad\quad {\mathrm{R}}(x,y)=
\frac{1}{4}{\mathrm{ad}}_{[x,y]}
$$
Hence every semi Riemannian bi-invariant metric is geodesically
complete, the geodesics through the unit element $\varepsilon$ of
$G$ are the 1-parameter subgroups of $G,$ and the bi-invariant
metric is flat if and only if the group is 2-step nilpotent.
\end{rmk}

\section{Jacobi fields on quadratic Lie groups, conjugate points}

\vspace{.1in} Every vector field $X$ on a Lie group defines a map
$\tilde{X}:G\to{\cal{G}}$ given by $\sigma\mapsto
({\mathrm{L}}_{\sigma^{-1}})_{*,\sigma}X_\sigma.$ Obviously, a
vector field is left invariant if and only if the associated map is
constant.

Given a curve $\sigma :[t_0,t_1]\to G,$ every vector field $Y$ on
$\sigma$ defines a curve in ${\cal{G}}:$
$$\tilde{Y}(t)=\left({\mathrm{L}}_{\sigma(t)^{-1}}\right)_{*,\sigma(t)}Y(t)
$$and conversely, every curve in ${\cal{G}}$ defined on $[0,1]$ determines a
vector field on $\sigma.$ We say that one is the {\bf{reflection}}
of the other and we write either $y^\sim=Y$ or $Y^\sim=y.$

Notice that $y(t)=(Y(t))^\sim$ is equivalent to
$y(t)_{\sigma(t)}^+=Y(t).$

The following proposition describes the Jacobi fields for a left
invariant semi Riemannian metric defined on a Lie group.

\begin{pr}\label{pr:jac}
Let $G$ be a Lie group, $\nabla$  a left invariant torsion free
connection on $G$  and let $\sigma :[0,1]\to G$ be a geodesic such
that $\sigma(0)=\varepsilon.$ Then the vector field on $\sigma,$
$t\mapsto Y(t)$ is a Jacobi vector field if and only if its
reflection $y:=\tilde{Y}$ satisfies the differential equation
\begin{equation}\label{eq:jac}
\ddot y +2 x\dot y= [y,x]x+ x[y,x]+[xx,y]
\end{equation}
where, as before,
$x(t)=\left({\mathrm{L}}_{\sigma^{-1}(t)}\right)_{*,\sigma(t)}\,\sigma'(t).$
\end{pr}

This proposition is a consequence of the following technical result.

\begin{lm}
With the notations introduced in the previous proposition, the first
and second covariant derivatives of the vector field $Y$ on $\sigma$
are given by
\begin{eqnarray*}
\frac{DY}{dt}&=&(y'+xy)^{\sim}\\
\frac{D^2Y}{dt^2}&=& (y''+2xy'+x'y+x(xy))^{\sim}
\end{eqnarray*}
\end{lm}

\vspace{.1in} \noindent {\bf{Proof.}} Let $G$ be a $n$ dimensional
Lie with Lie algebra ${\cal{G}}$ and  $\{e_i,\, 1\leq i\leq n\}$ be
a basis for the vector space underlying ${\cal{G}}.$ Expressing $Y$
by means of $e_i^+$ $(i=1,\cdots, n),$ as in proposition
\ref{pr:jac}, we get
\begin{eqnarray*}
\frac{DY}{dt}&=&\frac{D}{dt}\sum_{i=1}^n y_i(t)e_{i,\sigma(t)}^+\\
&=&\sum_{i=1}^n y_i^\prime(t)e_{i,\sigma(t)}^+ +\sum_{i=1}^n
y_i(t)\nabla_{\sigma'(t)}e_i^+\\
&=&(y'(t))^\sim +\sum_{i=1}^n y_i(t)\nabla_{\sigma'(t)}e_i^+ .
\end{eqnarray*}
Since $\nabla_{\sigma'(t)}e_i^+=(x(t)e_i)^+_{\sigma(t)},$
$$
\sum_{i=1}^n\, y_i(t)\nabla_{\sigma'(t)}e_i^+= \sum_{i=1}^n
y_i(t)(x(t)e_i)_{\sigma(t)}^+=
(x(t)y(t))_{\sigma(t)}^+=(x(t)y(t))^\sim,
$$
that is
$$
\frac{D}{dt}(\tilde{y})=(y'+xy)^{\sim}.
$$
Hence
\begin{eqnarray*}
\frac{D^2Y}{dt^2}&=&\frac{D}{dt}(y'+xy)^{\sim}\\
                 &=&\left( (y'+xy)'+ x(y'+xy)\right)^\sim\\
                 &=&\left( y'' +x'y+xy' +
xy'+x(xy)\right)^\sim\\
                 &=&\left( y'' +x'y+2xy' +x(xy)\right)^\sim
\end{eqnarray*}
\qed

\vspace{.1in} \noindent {\bf{Proof of Proposition \ref{pr:jac}.}}
Recall that
$$
{\mathrm{R}}_{x^+y^+}z^+=\left([x,y]z-x(yz)+y(xz)\right)^+.
$$
A straightforward calculation shows that,
$$
{\mathrm{R}}_{Y\sigma'(t)}\sigma'(t)= ([y,x]z-y(xz)+x(yz))^\sim.
$$
Thus $Y=y^\sim$ is a Jacobi vector field if and only if
\begin{equation}\label{eq:jac2}
\ddot y +\dot{x}y+2x\dot{y}+x(xy) =[y,x] z -
 y  (x z )+  x  (y z ).
\end{equation}

Since $\sigma$ is a geodesic $\dot x=-  x x, $ and since the
connection is torsion free, $  x y -  y x =[x,y],$ then equation
(\ref{eq:jac2})  becomes
\begin{eqnarray*}
\ddot y +2  x \dot y  &=&
  [y,x] x +    (x x) y -  y  (x x)  +  x [y,x] \\
&=&  [y,x] x +  x [y,x] +  [(x x), y ]
\end{eqnarray*}
\qed

\begin{rmk}
If $x_0\ne 0$ is a solution of $x x=0,$ then the geodesic $\sigma$
through $\varepsilon$ with velocity $x_0$ is the one-parameter
subgroup of $G$ with infinitesimal generator $x_0$ and the
reflections of Jacobi fields on $\sigma$ are the solutions of the
equation
\[
\ddot y +2  x_0 \dot y =  [y,x_0] x_0 +  x_0 [y,x_0]
\]
\end{rmk}

\begin{cor}\label{cor:jac}
If $\nabla$ is the Levi-Civita connection defined by a bi-invariant
semi Riemannian metric, then a vector field $Y$ on a  geodesic
$\sigma :[0,1]\to G$ with $\sigma(0)=\varepsilon,$ is a Jacobi
vector field is and only if its reflection curve $y=Y^\sim$ is a
solution of the differential equation
$$
\ddot y= [\dot y,x_0],
$$
where $x_0$ is the initial velocity of the geodesic.
\end{cor}

\vspace{.1in} \noindent The proof follows immediately from the fact
that the Levi-Civita product of a bi-invariant metric is given by $
x y =(1/2)[x,y].$

\begin{cor}\label{cor:jacnilp}
Let $\nabla$ be the Levi-Civita connection defined by a bi-invariant
semi Riemannian metric on a  nilpotent Lie group. Then the
reflection of Jacobi fields along a  geodesic $\sigma :[0,1]\to G$
are polynomial.
\end{cor}

\vspace{.1in} \noindent{\bf{Proof.}} Let $Y$ be a Jacobi field along a
geodesic $\sigma,$ and $y$ its reflection. Then by the previous corollary
 $ y^{(2)}(t)= [y'(t),x_0].$ Hence $y^{(k+1)}(t)=[y^{(k)}(t),x_0],$ and
 $y^{(m+1)}\equiv 0,$ where $m$ is such that ${\cal{G}}^{(m)}=0.$
\qed

\subsection*{Jacobi fields on the oscillator groups}

Consider the quadratic Lie algebra  $(\mathbb{R}^{2n},k_0)$, where
$k_0$ is an Euclidean inner product. Let
$\lambda:=(\lambda_1,\ldots,\lambda_n)$ where each $\lambda_i$ is a
positive real number and $\theta$ the antisymmetric isomorphism of
$(\mathbb{R}^{2n},k_0)$ given by the matrix (with respect to an
orthonormal basis)
${\cal{B}}=\{e_1,\ldots,e_n,\check{e}_1,\ldots,\check{e}_n)$
$$
M_{\cal{B}}\theta=\left(
\begin{array}{cccc}
0 & -{\mathrm{diag}}(\lambda_1,\ldots,\lambda_n)\\
{\mathrm{diag}}(\lambda_1,\ldots,\lambda_n) & 0
\end{array}
\right)
$$
where ${\mathrm{diag}}(\lambda_1,\ldots,\lambda_n) $ stands for the
diagonal matrix with $\lambda_1,\ldots,\lambda_n$ on the main
diagonal.

Obviously $\theta$ defines a representation of the Lie algebra
$\mathbb{R}$ by endomorphisms of ${\cal{G}},$ noted also by
$\theta.$

Let ${\cal{G}}(\lambda)$ be Lie algebra obtained by a process of
double extension (see \cite{kn:MRMM}) of $(\mathbb{R}^{2n},k_0)$ by
$\mathbb{R}$ by means of $\theta.$ This means that the algebra is
obtained by a central extension
$\mathbb{R}\,e_0\times_\omega\mathbb{R}^{2n}$ of $\mathbb{R}^{2n}$
by $\mathbb{R}$ by means of the scalar 2-cocycle $\omega (x,y):=
k_0(\theta (x),y),$ then by semi-direct product of $\mathbb{R}\,
e_{-1}$ by $\mathbb{R}\,e_0\times_\omega\mathbb{R}^{2n}$ where the
action is given by
$$
[e_{-1},e_0]=0,\;\; [e_{-1},x]=\theta x, \;{\mbox{for}}\;
x\in\mathbb{R}^{2n}.
$$
This algebra has a quadratic structure  $k$ that extends $k_0$ and
is given in the Minkowski plane, $V={\mathrm{Span}}\{e_{-1},e_0\},$
by
$$
\begin{array}{rc}
k(e_0,e_0)=k(e_{-1},e_{-1})=0 , & k(e_{-1},e_0)=1
\end{array}
$$
and is orthogonal to $\mathbb{R}^{2n}.$

Since the algebra ${\cal{G}}(\lambda)$ is solvable, the connected
and simply-connected Lie group with Lie algebra
${\cal{G}}(\lambda_1,\ldots,\lambda_n)$ can be identified as a
manifold to $\mathbb{R}^{2n+2}\equiv
\mathbb{R}\times\mathbb{C}^n\times\mathbb{R}$ with product
$$
(s,z_1,\ldots,z_n,t)\cdot(s',z'_1,\ldots,z'_n,t')=
$$
$$
(s+s'+\frac{1}{2}\sum_{j=1}^n{\mathrm{Im}}\bar{z_j}\exp(i\,
t\lambda_j)z'_j, z_1+\exp(i\, t\lambda_1)z'_1, \ldots, z_n+\exp(i\,
t\lambda_n)z'_n,t+t')
$$

\begin{df}
The groups ${\mathrm{G}}(\lambda)$ are called {\bf{Oscillator
groups}} and the corresponding Lie algebras ${\cal{G}}(\lambda)$ are
called {\bf{ Oscillator algebras}}.
\end{df}

The equation that defines the reflection of a Jacobi field  in the
oscillator algebra is given by
$$
\begin{array}{lcl}
y_{-1}^{\prime\prime}&=&0\\
y_0^{\prime\prime}&=&
\check{x}_1 y_1^\prime+ \cdots +\check{x}_ny_n^\prime- x_1\check{y}^\prime_1- \cdots -
x_n\check{y}^\prime_n\\
y_1^{\prime\prime}&=&-\lambda_1\check{x}_1y_{-1}^\prime+x_{-1}\lambda_1\check{y}^\prime_1\\
 & \ldots & \\
y_n^{\prime\prime}&=&-\lambda_n\check{x}_n y_{-1}^\prime+x_{-1}\lambda_n\check{y}^\prime_n\\
\check{y}_1^{\prime\prime}&=&\phantom{-}\lambda_1x_1y_{-1}^\prime-x_{-1}\lambda_1y_1^\prime \\
&\ldots &\\
\check{y}_n^{\prime\prime}&=&\phantom{-}\lambda_nx_ny_{-1}^\prime-x_{-1}\lambda_ny_n^\prime,
\end{array}
$$
where $ x(0)=\sum_{i=-1}^n\, x_i\,e_i + \sum_{i=1}^n\,
\check{x}_i\,\check{e}_i. $ In order to find the conjugate points to
$\varepsilon,$ it is also necessary that  $y(0)=0$ and $y(t_1)=0$
for some $t_1\neq 0.$ This implies that $y_{-1}\equiv 0$ and the
system is equivalent to the system
$$
\begin{array}{lcl}y_0^{\prime\prime}&=&\check{x}_1y_1^\prime+
\cdots +\check{x}_ny_n^\prime - x_1\check{y}^\prime_1- \cdots -
x_n\check{y}^\prime_n\\
y_j^{\prime\prime}&=& x_{-1}\lambda_j\check{y}^\prime_j\\
\check{y}_j^{\prime\prime}&=&-x_{-1}\lambda_jy_j^\prime \\
\end{array},
$$
$1\leq j\leq n.$ When $x_{-1}\neq 0,$ the solutions are
$$
\begin{array}{lcl}
y_j(t)&=&\displaystyle\frac{r_j}{x_{-1}\lambda_j}\sin(x_{-1}\lambda_j\,t)\\
&&\\
\check{y}_j(t)&=&\displaystyle\frac{r_j}{x_{-1}\lambda_j}
\left(1-\cos(x_{-1}\lambda_j\,t)\right),
\end{array}
$$
since $y(0)=0.$ In order to have $y(t_1)=0,$ it is necessary that
$$
x_{-1}\lambda_j\,t_1=2\pi k,\;\; k\in\mathbb{Z},\;\;{\mbox{or}}\;\;
r_j=0
$$
for all $1\leq j\leq n.$ Hence, letting  $r_j=0$ for $j\neq i,$
$y_0$ must be a solution to the equation
\begin{eqnarray*}
y_0^{\prime\prime}&=&
\check{x}_iy_i^\prime- x_i\check{y}^\prime_i\\
&=&r_i \left(\check{x}_i\cos(x_{-1}\lambda_i\,t)-
x_i\sin(x_{-1}\lambda_i\,t)\right)
\end{eqnarray*}
which implies that
\[
y_0'(t)=c+\frac{r_i}{x_{-1}\lambda_i}
\left(x_i\cos(x_{-1}\lambda_i\,t)+\check{x}_i
\sin(x_{-1}\lambda_i\,t)\right)\] and thus \[y_0(t)=d+ct
+\frac{r_i}{(x_{-1}\lambda_i)^2}
\left(x_i\sin(x_{-1}\lambda_i\,t)-\check{x}_i
\cos(x_{-1}\lambda_i\,t)\right).
\]
Since $y_0(0)=y_0(t_1)=0,$
$$
y_0(t)=\frac{r_i}{(x_{-1}\lambda_i)^2}
\left(x_i\sin(x_{-1}\lambda_i\,t)+\check{x}_i(1- \cos(x_{-1}\lambda_i\,t))\right)\\
$$
and the vector field $Y=y^\sim,$ where
$$
y(t)=y_0(t)e_0+\frac{r_i}{x_{-1}\lambda_i}\sin(x_{-1}\lambda_i t)e_i
+
\frac{r_i}{x_{-1}\lambda_i}\left(1-\cos(x_{-1}\lambda_i\,t)\right)\check{e}_i,
$$
is a Jacobi field along the geodesic $t\mapsto\exp(tx(0)).$ The
points $\exp(\displaystyle\frac{\pi k}{x_{-1}\lambda_i}x(0)),$
$k\in\mathbb{Z}, 1\leq i\leq n,$ are conjugate to $\varepsilon.$

Notice that when $x_{-1}=0$ then a Jacobi field along the geodesic
$t\mapsto\exp(tx(0)),$ with $y(0)=0$ vanishes everywhere.\qed

In what follows, $G$ is a connected Lie group.

\section{On flat left invariant semi Riemannian metrics on quadratic Lie groups}

Let $G$ be a quadratic Lie group and $\langle ,\rangle$ a flat left
invariant  semi Riemannian metric. Notice that a  quadratic Lie
group is unimodular because ${\mathrm{ad}}_x$ is $k$ skew-symmetric,
for all $x\in{\cal{G}},$ hence, by Theorem \ref{theorem:uni}, every
flat left invariant semi Riemannian metric on a quadratic Lie group
is complete, and we have a companion theorem of Theorem
\ref{theorem:ACH},

\begin{theorem}\label{theorem:ACHQ}
Let $G$ be a connected quadratic Lie group, $\langle\, ,\,\rangle$ a
left invariant semi Riemannian metric on $G.$ Then the following
assertions are equivalent

\noindent i) $\langle\, ,\,\rangle$ is flat.

\noindent ii) the exponential map relative to $\langle\, ,\,\rangle$
is a local isometry.

\noindent iii)  $\langle\, ,\,\rangle$ is flat and complete.

In any case,  $G$ is solvable, and $(G,\langle\,,\rangle)$ has no
conjugate points.

Moreover $G$ viewed as a group of transformations of ${\cal{G}}$ has
non trivial central 1-parameter subgroups of translations.

\end{theorem}

\vspace{.1in}\noindent{\bf{Proof.}} A quadratic Lie group is
unimodular, hence the first part of the theorem follows from Theorem
\ref{theorem:ACH}.

As for the existence of non-trivial 1-parameter subgroups of
translations, since the metric is complete, the Levi-Civita product
has no non trivial idempotents, hence there is an element
$x_0\in{\cal{G}},$ $x_0\ne 0$ such that $x_0x_0=0.$ The 1-parameter
subgroup with infinitesimal generator $x_0$ is a geodesic for the
metric. \qed

\begin{df}
The {\bf{index}} of a non degenerate quadratic form $q$ defined on a
real vector space $V$ is the maximal dimension of a $q$ totally
isotropic subspace of $V.$
\end{df}

The principal result of this section is

\begin{theorem}\label{theorem:ind>2}
If a connected non Abelian  quadratic Lie group $(G,k)$ admits a
flat left invariant semi Riemannian metric $\langle,\rangle$, then
$\langle,\rangle$ is geodesically complete, $G$ is solvable, the
index of $\langle,\rangle$ is $>1$ and the universal covering of $G$
viewed as a group of affine transformations contains central
1-parameter groups of translations.
\end{theorem}

The two first assertions of the theorem follow from Theorem
\ref{theorem:ACH}. The other assertions will follow from a series of
lemmas and propositions.

\begin{lm}\label{lm:center}
The center of a quadratic Lie  group with a flat left invariant semi
Riemannian metric is non trivial.
\end{lm}

\noindent{\bf{Proof.}} By Theorem \ref{theorem:ACH} the Lie group is
solvable. The result follows from the observation that
$\mathcal{Z}({\cal{G}})^{\perp_k}=[{\cal{G}},{\cal{G}}].$ \qed

\begin{pr}\label{pr:lf} Let $(G,k)$ be quadratic Lie group with a
flat left invariant semi Riemannian metric. Then then for all
$e\in{\cal{Z}}\, ({\cal{G}}),$ $\nabla_e^2=0.$ If the metric is
either Riemannian or Lorentzian, then for all $e\in{\cal{Z}}\,
({\cal{G}}),$ $\nabla_e=0.$ In this case, if
$u\in{\mathrm{Gl}}({\cal{G}})$ is the $k$ symmetric isomorphism
induced by the semi Riemannian  metric, then ${\cal{Z}}\,
({\cal{G}})$ is invariant by $u.$
\end{pr}

\noindent {\bf{Proof.}} Consider the Levi-Civita product induced by
the semi Riemannian metric. It is immediate from the Koszul formula
that $e e' =L_ee'= 0$ for $e,e'\in{\cal{Z}}({\cal{G}}).$ The
following string of equalities
\begin{eqnarray*}\label{eq:imL}
 \langle L_eL_e x,y\rangle =-\langle L_e x ,L_e y\rangle &=&
-\langle L_e x ,L_y e \rangle =\langle L_yL_e x ,e \rangle
\\&=&\langle L_e L_y x,e\rangle =-\langle L_y x,L_e e\rangle =0,
\end{eqnarray*}
implies the first assertion and it also implies that the subspace
${\mathrm{Im}}\,(L_e)$ is totally isotropic. If the metric is either
Riemannian or Lorentzian ${\mathrm{dim}}\,{\mathrm{Im}}\,(L_e)\le
1.$ Suppose that there exists $e\in{\cal{Z}}({\cal{G}})$ such that
$L_e\ne 0,$ and let $x$ such that $L_ex\ne 0.$ Then for every
$y\in{\cal{G}},$ there is a $\lambda\in\mathbb{R}$ such that
$L_ey=\lambda L_ex.$ Notice that $\langle x,L_ex\rangle =0$ because
$L_e$ is $\langle\, ,\rangle$ skew symmetric. Then
$$
0=\langle L_ey,x\rangle +\langle y,L_ex\rangle = \lambda\langle
L_ex,x\rangle +\langle y,L_ex\rangle = \langle y,L_ex\rangle ,$$
thus, for all $y \in {\cal{G}},$
$$
\langle y,L_ex\rangle =0.
$$
This equality implies that $L_ex=0,$ because the semi Riemannian
metric is non degenerate, contrary to the assumption.

For the second part of the assertion, the Koszul formula
$$
\langle  x y ,z\rangle =
\frac{1}{2}\left(\langle[x,y],z\rangle -\langle[y,z],x\rangle +\langle[z,x],y\rangle \right)
$$
for $x=e,$  $e\in{\cal{Z}}({\cal{G}}),$  reduces to
$$
\langle L_ey,z\rangle =-\frac{1}{2}\langle[y,z],e\rangle
$$
for all $y,z\in{\cal{G}}.$ Let $u$ be the $k$ symmetric isomorphism
of the underlying vector space of ${\cal{G}}$ induced by the semi
Riemannian metric $\langle\, ,\rangle .$ Then, using that $L_e=0,$
$$
0=\langle L_e(y),z\rangle =-\frac{1}{2}\langle[y,z],e\rangle =
-\frac{1}{2}k([y,z],u(e))=-\frac{1}{2}k([u(e),y],z).
$$
The former equality implies that
$$
0=[u(e),y],
$$
because $k$ is non degenerate. Thus, for a flat left invariant semi
Riemannian metric of index $<2,$ $u(e)\in{\cal{Z}}({\cal{G}}).$

\qed

\begin{cor}
Under the hypothesis of proposition \ref{pr:lf} if $G$ is viewed as
a group of affine transformations of ${\cal{G}},$  it has
one-parameter subgroups of translations.
\end{cor}

As announced in Section 2, the existence of flat left invariant
Riemannian metrics on quadratic Lie groups imposes severe
restrictions on the group:

\begin{pr}\label{pr:fq}
A quadratic Lie group with a flat left invariant Riemannian metric
is Abelian.
\end{pr}

\noindent {\bf{Proof.}} By Lemma \ref{lm:center},
${\cal{Z}}({\cal{G}})\ne 0.$ Consider the map
$L:{\cal{G}}\to{\mathit{gl}}({\cal{G}})$ defined by $x\mapsto L_x.$
The fact that $[x,y]=L_xy-L_yx$ implies that ${\mathrm{ker}}\,(L)$
is an Abelian ideal of ${\cal{G}}.$ By Proposition \ref{pr:lf},
${\cal{Z}}({\cal{G}})\subset {\mathrm{ker}}\,(L).$ Since the metric
is flat, by Theorem 1.5 in \cite{kn:Mi}
$$
{\cal{G}}={\mathrm{ker}}\,(L)\oplus {\cal{H}}
$$
where ${\cal{H}}:= {\mathrm{ker}}\,(L)^\perp$ and ${\cal{H}}$ acts
on ${\mathrm{Ker}}(L)$ by adjoints. Hence
$[{\cal{G}},{\cal{G}}]\subset {\mathrm{ker}}\,(L).$ Since
$[{\cal{G}},{\cal{G}}]$ is orthogonal to ${\cal{Z}}({\cal{G}})$
relative to $\langle\,,\,\rangle $ (if $e\in{\cal{Z}}({\cal{G}}),$
and $x,y\in{\cal{G}}$ then $\langle[x,y],e\rangle
=k([x,y],u(e))=k(x,[y,u(e)])=0,$ because ${\cal{Z}}({\cal{G}})$ is
invariant by $u.)$

Since ${\cal{G}}$ is quadratic,
$${\mathrm{dim}}{\cal{Z}}({\cal{G}})+
{\mathrm{dim}}[{\cal{G}},{\cal{G}}]={\mathrm{dim}}{\cal{G}}.$$

Hence ${\cal{H}}=(0).$ \qed

\begin{pr}
A flat left invariant  semi Riemannian metric on a non Abelian
quadratic Lie group has index $\ge 2.$
\end{pr}

\noindent {\bf{Proof.}} Let $({\cal{G}},k)$ a non Abelian  quadratic
Lie algebra. Then, by Proposition \ref{pr:fq}, ${\cal{G}}$ has no
flat left invariant Riemannian metric. Suppose that $u$ is a $k$
symmetric isomorphism of ${\cal{G}}$ such that $k_u:=\langle\,
,\,\rangle $ is Lorentzian and flat.  We proceed by induction on the
dimension of ${\cal{G}}.$ A non Abelian quadratic algebra of
dimension $\le 3$ has no flat left invariant semi Riemannian metric.
Since ${\cal{Z}}({\cal{G}})\ne (0),$ suppose first that
${\cal{Z}}({\cal{G}})$ is non degenerate for $k.$ Then
$$
{\cal{G}}={\cal{Z}}({\cal{G}})\oplus_{\perp_k} [{\cal{G}},{\cal{G}}]
$$
because, as noted before
$[{\cal{G}},{\cal{G}}]={\cal{Z}}({\cal{G}})^{\perp_k}.$ Furthermore,
$[{\cal{G}},{\cal{G}}]$ is a non Abelian quadratic Lie algebra.
Since $u$ is $k$ symmetric and leaves invariant
${\cal{Z}}({\cal{G}}),$ it leaves also invariant
$[{\cal{G}},{\cal{G}}],$ and $\langle\, ,\,\rangle$ is non
degenerate both on ${\cal{Z}}({\cal{G}})$ and
$[{\cal{G}},{\cal{G}}].$ The restriction of $\langle\, ,\,\rangle$
to $[{\cal{G}},{\cal{G}}]$ is either Riemannian or Lorentzian.
Because of Proposition \ref{pr:fq} and by the induction hypothesis
it cannot be either one. Hence ${\cal{Z}}({\cal{G}})$ is degenerate
for $k.$ Let $e_0\in{\cal{Z}}({\cal{G}})$ such that $k(e_0,e)=0$ for
all $e\in{\cal{Z}}({\cal{G}}) ,$ that is
$e_0\in{\cal{Z}}({\cal{G}})\cap [{\cal{G}},{\cal{G}}].$ By
Proposition \ref{pr:lf}, $u(e_0)\in {\cal{Z}}({\cal{G}}),$ hence
$\langle e_0 ,[x,y]\rangle = k(u(e_0),[x,y])= k([u(e_0),x],y)=0$ for
all $x,y\in{\cal{G}}.$ This implies that
${\mathrm{Span}}\,\{e_0,u(e_0)\}$ is $\langle\, ,\,\rangle$ totally
isotropic. The semi Riemannian metric being Lorentzian, the vector
space is of dimension 1, hence $u(e_0)=\nu e_0,$ for some $\nu\ne
0.$

Since $[{\cal{G}},{\cal{G}}]$ is invariant by $u,$ the algebra
$[{\cal{G}},{\cal{G}}]/\mathbb{R} e_0$ is a quadratic Lie algebra
with a flat left invariant Riemannian metric. By Proposition
\ref{pr:fq} it is Abelian, that is, $\forall\,
x,y\in[{\cal{G}},{\cal{G}}],\; [x,y]\in\mathbb{R} e_0.$ The
isomorphism $u$ induces an isomorphism $\tilde{u}$ on
$[{\cal{G}},{\cal{G}}]/\mathbb{R} e_0,$ and the reduced metric
induced by $\langle\, ,\,\rangle$ on
$[{\cal{G}},{\cal{G}}]/\mathbb{R} e_0,$ is positive definite. Then
there is an  orthonormal basis, $\{\tilde{E}_i\},$ relative to the
reduced metric, that diagonalizes $\tilde{u}:$
$\tilde{u}(\tilde{E}_i)=\nu_i\tilde{E}_i.$

Let $E_i$ an element in $[{\cal{G}},{\cal{G}}]$ that projects onto
$\tilde{E}_i$ in $[{\cal{G}},{\cal{G}}]/\mathbb{R} e_0.$

Then $\{e_0, E_1,\cdots, E_{n}\}$ is a basis of
$[{\cal{G}},{\cal{G}}]$ and
$$
\langle E_i , E_j \rangle=\delta_{ij}\quad\quad k(e_0,E_i)=0\quad
i,j\ge 1.
$$
Let
$$
u(E_i)=\nu_iE_i+\mu_ie_0.
$$
($\mu_i\in\mathbb{R}\,,i\ge1),$ then
$$
k(E_i,E_j)=\frac{\delta_{ij}}{\nu_i}.
$$

The same construction can be done for ${\cal{Z}}({\cal{G}})\slash
\mathbb{R} e_0.$ The isomorphism $u$ induces an isomorphism
$\tilde{u}$ on ${\cal{Z}}({\cal{G}}) /\mathbb{R} e_0,$ and the
reduced metric induced by $k_u$ on ${\cal{Z}}({\cal{G}})/\mathbb{R}
e_0,$ is positive definite. Then there is an orthonormal basis
relative to the induced metric $\{\tilde{F}_i\},$ that diagonalizes
the induced isomorphism: $\tilde{u}(\tilde{F}_i)=\nu'_i\tilde{F}_i.$

Let $e_{-1}$ any vector in ${\cal{G}}$ not in ${\cal{Z}}({\cal{G}})
+ [{\cal{G}},{\cal{G}}].$ Then the vector space underlying the Lie
algebra ${\cal{G}}$ decomposes as:
$$
{\cal{G}}= {\mathrm{Span}}\,\{F_i\,:\, 1\le i\le m\}\oplus
{\mathrm{Span}}\,\{E_i\,:\, 1\le i\le n\} \oplus
{\mathrm{Span}}\,\{e_{-1},e_0\}.
$$
The vector $e_{-1}$ can be chosen to satisfy $\langle e_{-1}, e_0\rangle =1$ and
$$
e_{-1}\perp_{\langle \,\rangle}{\mathrm{Span}}\,\{F_i\,:\, 1\le i\le
m\} \quad e_{-1}\perp_{\langle \,\rangle}{\mathrm{Span}}\,\{E_i\,:\,
1\le i\le n\}.
$$

Since $u$ is $k$-symmetric,
$$
u(e_{-1})=\nu e_{-1}+\rho e_0+\sum \mu_iE_i + \sum \mu'_iF_i
$$
where $\rho=\langle e_{-1},e_{-1}\rangle .$

Denote by $[E_i,E_j]=\rho_{ij}e_0,$ (for $i,j\ge 1)$ and notice that
$\rho_{ij}=-\rho_{ji}.$ In order to calculate
$$
\langle{\mathrm{R}}_{e_{-1}E_i}E_i,e_{-1}\rangle
$$
some remarks are needed:
\begin{enumerate}
\item
$\forall i\ge 1,\quad L_{E_i}E_i=0, $ because
$$\langle L_{E_i}E_i,x\rangle =\langle[x,E_i],E_i\rangle =k([x,E_i],\nu_iE_i+\mu_i e_0\rangle =0,
\forall x\in{\cal{G}}.$$
\item
$\forall x,y\in{\cal{G}},\quad \langle L_xy,y\rangle =0.$
\item
$\forall x,y\in{\cal{G}},\quad  \langle L_xy,e_0\rangle =0.$
\item
$[e_{-1},E_i]=\sum_j\rho_{ij}E_j$ because
$[e_{-1},E_i]=ae_0+\sum_j\alpha_jE_j$ and
$k([e_{-1},E_i],E_j)=k(e_{-1},[E_i,E_j]) =\rho_{ij},$
$a=k[e_{-1},E_i],e_{-1})=0.$
\item $\langle L_{E_j}e_{-1},E_i\rangle =$
$$\begin{array}{l}
=(1/2)\{\langle[E_j,e_{-1}],E_i\rangle -\langle[e_{-1},E_i],E_j\rangle +\langle[E_i,E_j],e_{-1}\rangle \\
=(1/2)(-\rho_{ji}\nu_i-\rho_{ij}\nu_j+\rho_{ij}\nu)\\
=(\rho_{ij}/2)(\nu_i-\nu_j+\nu)
\end{array}$$
\end{enumerate}

After a  rather cumbersome, yet straightforward calculation, we get
$$\langle{\mathrm{R}}_{e_{-1}E_i}E_i,e_{-1}\rangle =$$
$$
\begin{array}{l}
\langle L_{[e_{-1},E_i]}E_i,e_{-1}\rangle
-\langle L_{e_{-1}}L_{E_i}E_i-L_{E_i}L_{e_{-1}}E_i,e_{-1}\rangle =\\
-\left(\sum_j\rho_{ij}\langle E_i,L_{E_j}e_{-1}\rangle \right)-\langle L_{e_{-1}}E_i,L_{E_i}e_{-1}\rangle =\\
-(1/2)\sum_j\rho_{ij}^2(\nu_i-\nu_j+\nu)-\langle L_{E_i}e_{-1},L_{E_i}e_{-1}+[e_{-1},E_i]\rangle .
\end{array}
$$
Since
$$
\langle L_{E_i}e_{-1},L_{E_i}e_{-1}\rangle =
\dps\sum_j\dps\frac{\rho_{ji}^2}{4\nu_j}(\nu_j-\nu_i+\nu)^2
$$
and
$$
\langle L_{E_i}e_{-1},[e_{-1},E_i]\rangle =\dps\sum_j\langle
L_{E_i}e_{-1},\rho_{ij}E_j\rangle =
\dps\sum_j\dps\frac{\rho_{ij}^2}{2}(\nu_j-\nu_i+\nu),
$$
we have
$$
\langle{\mathrm{R}}_{e_{-1}E_i}E_i,e_{-1}\rangle
=\sum_j\rho^2_{ij}(\nu_j-\nu_i)
-\sum_j\frac{\rho^2_{ij}}{4\nu_j}(\nu_j-\nu_i+\nu)^2
$$
and
$$
\sum_i\langle{\mathrm{R}}_{e_{-1}E_i}E_i,e_{-1}\rangle =
-\sum_{i,j}\frac{\rho^2_{ij}}{4\nu_j}(\nu_j-\nu_i+\nu)^2.
$$
Recall that the metric is Lorentzian and flat, then $\nu_j > 0$ and
$$
0= \sum_{i,j}\frac{\rho^2_{ij}}{4\nu_j}(\nu_j-\nu_i+\nu)^2,
$$
implies that, for all $i,j\ge 1,$
$$
0= \rho^2_{ij}(\nu_j-\nu_i+\nu)^2.
$$
Notice that $\rho_{ij}\ne 0,$ for some $i,j$ for otherwise
$[{\cal{G}},{\cal{G}}]$ would be Abelian. Then
$$
\begin{array}{ccc}
\nu_j-\nu_i+\nu&=&0\\
\nu_i-\nu_j+\nu&=&0
\end{array}
$$
and $\nu=0,$ which contradicts the hypothesis. \qed

\vspace{.1in}\noindent{\bf{Proof of Theorem \ref{theorem:ind>2}.}}
Let $G$ be a quadratic non Abelian Lie group. By Theorem
\ref{theorem:ACHQ} every flat left invariant semi Riemannian metric
on $G$ is geodesically complete, and by Theorem \ref{theorem:ACH} ,
$G$ is solvable. By Proposition 5, the geodesic through the unit of
$G$ with velocity $e_0$ is the 1-parameter subgroup of $G$ with
infinitesimal generator $e_0,$ because $e_0e_0=0.$ Finally,
Proposition 7 states that the index of the metric is $\ge 2.$

The following corollaries are also consequences of Theorem
\ref{theorem:ind>2}:

\begin{cor}
Every left invariant semi Riemannian metric on a quadratic Lie group
with non trivial Levi component (in particular when the group is reductive)  is non flat.
\end{cor}

\begin{df}
A quadratic Lie group $(G,k)$ is called {\bf{undecomposable}} if it
has no non trivial normal Lie subgroups $N$ such that the
restriction of the bi-invariant metric $k$ to $N$ is non degenerate.
At the algebra level, this means that every ideal of ${\cal{G}}$ is
$k$ degenerate.
\end{df}

\begin{cor}
Let $(G,k)$ be an undecomposable quadratic Lie group that admits a flat
left invariant semi Riemannian metric. Then the index of $k$ is $\ge 2.$
\end{cor}

\begin{theorem}\label{theorem:dim4} Every undecomposable quadratic Lie group of
dimension 4 has affine left invariant structures and no flat left
invariant semi Riemannian metrics.
\end{theorem}

\vspace{.1in}\noindent{\bf{Proof.}} There are two undecomposable
quadratic connected Lie groups. The first one is the oscillator
group of dimension 4 that was treated in \cite{kn:bm2}. The Lie
algebra of the second one is  obtained as follows.

Let ${\cal{G}}$ be the Abelian Lie algebra obtained by quadratic
double extension from the Minkowski plane
$\mathrm{Span}\,\{e_1,e_2\}$ (viewed as an Abelian Lie algebra) by a
central line $\mathbb{R}\,e_0.$ This is a four-dimensional quadratic
Lie algebra with an $\mathrm{ad}$ antisymmetric scalar product of
index 2. It has a basis $e_{-1},e_0,e_1,e_2$ with bracket
$$
[e_1,e_2]=-e_0,\quad [e_{-1},e_1]=e_2,\quad [e_{-1},e_2]=e_1
$$
and quadratic structure
$$
k(e_{-1},e_0)=k(e_1,e_1)=-k(e_2,e_2)=1,
$$
the non stated products are either given by antisymmetry/symmetry or
are 0. Note that ${\cal{Z}}({\cal{G}})=\mathbb{R} e_0.$

We claim that ${\cal{G}}$ has no left invariant flat semi Riemannian
metric. Suppose on the contrary that there exists a $k$ symmetric
isomorphism $u$ of the vector space underlying ${\cal{G}}$ such that
the metric
$$
\langle x,y\rangle:=  k(x ,u(y))
$$
is flat. By Proposition \ref{pr:lf} , we have that
$$
{\mathrm{im}}\,L_{e_0}\subset{\mathrm{ker}}\,L_{e_0}.
$$
Recall that, for a left invariant semi Riemannian metric, the
Levi-Civita product is given by:
$$
2  x y =[x,y]+u^{-1}([x,u(y)]+[y,u(x)]).
$$
Then
$$
2L_{e_0}=-u^{-1}\circ {\mathrm{ad}}_{u(e_0)}
$$
hence
\begin{equation}\label{eq:kim}
{\mathrm{im}}\,{\mathrm{ad}}_{u(e_0)}=u({\mathrm{im}}\,L_{e_0})\subset
u({\mathrm{ker}}\,L_{e_0})=u({\mathrm{ker}}\,{\mathrm{ad}}_{u(e_0)}).
\end{equation}
This equation implies that
$$
{\mathrm{rank}}\,{\mathrm{ad}}_{u(e_0)}\le 2.
$$
Let $u(e_0)=x_{-1}e_{-1}+x_{0}e_{0}+x_{1}e_{1}+x_{2}e_{2}.$

If $x_{-1}\ne 0,$ then $e_0$ and $u(e_0)$ are linearly independent
and
\begin{equation}\label{eq:laotra2}
{\mathrm{im}}\,{\mathrm{ad}}_{u(e_0)}=u({\mathrm{ker}}\,{\mathrm{ad}}_{u(e_0)})
=\mathrm{Span}\{u(e_0),u^2(e_0)\}.
\end{equation}
It is easy to show that
$${\mathrm{im}}\,{\mathrm{ad}}_{u(e_0)}=
\mathrm{Span}\{u(e_1)=
x_{2}e_{0}+x_{-1}e_{2},u(e_2)=-x_{1}e_{0}+x_{-1}e_{1}\}.
$$
By Equation \ref{eq:laotra2}, there exist $A_1,B_1,A_2,B_2$ such
that
\begin{eqnarray*}
x_{2}e_{0}+x_{-1}e_{2}&=& A_1u(e_0)+B_1u^2(e_0)\\
-x_{1}e_{0}+x_{-1}e_{1}&=& A_2u(e_0)+B_2u^2(e_0)
\end{eqnarray*}
Then
$$
\begin{array}{ccccc}
0&=&k(x_{2}e_{0}+x_{-1}e_{2},e_0)&=& k(A_1u(e_0)+B_1u^2(e_0),e_0)
\\
0&=&k(-x_{1}e_{0}+x_{-1}e_{1},e_0)&=&k( A_2u(e_0)+B_2u^2(e_0),e_0).
\end{array}
$$
and
$$
\begin{array}{ccc}
0&=&A_1k(u(e_0),e_0)+B_1k(u^2(e_0),e_0)\\
0&=& A_2k(u(e_0),e_0)+B_2k(u^2(e_0),e_0).
\end{array}
$$
 As a consequence, $A_1=\lambda A_2,$   $B_1=\lambda B_2,$ and
$$
x_{2}e_{0}+x_{-1}e_{2}= A_1u(e_0)+B_1u^2(e_0)=\lambda
(A_2u(e_0)+B_2u^2(e_0))=\lambda (-x_{1}e_{0}+x_{-1}e_{1}).
$$
This equality contradicts the fact that $e_0, e_1, e_2$ are linearly
independent.

If $x_{-1}=0,$ and $x_1^2+x_2^2\ne 0,$ then $u(e_0)=x_0
e_0+x_1e_1+x_2e_2, $ and
$\mathrm{dim}\,({\mathrm{ker}}\,{\mathrm{ad}}_{u(e_0)})=2.$ Hence
$\mathrm{dim}\,({\mathrm{im}}\,{\mathrm{ad}}_{u(e_0)})=2,$ and
$$
{\mathrm{im}}\,({\mathrm{ad}}_{u(e_0)})=u({\mathrm{ker}}\,{\mathrm{ad}}_{u(e_0)}).
$$
Moreover
${\mathrm{ker}}\,({\mathrm{ad}}_{u(e_0)})={\mathrm{Span}}\,\{e_0,u(e_0)\}$
and $e_0\in{\mathrm{im}}\,{\mathrm{ad}}_{u(e_0)}.$

By (\ref{eq:kim}), $e_0=\lambda u(e_0)+\mu u^2(e_0)$ with $\mu\ne
0,$ because we are supposing that $u(e_0)\notin\mathbb{R} e_0.$ Then
$$
0=k(e_0,e_0)=k(\lambda u(e_0)+\mu u^2(e_0),e_0)=\lambda
k(u(e_0),e_0)+ \mu k( u^2(e_0),e_0)
$$
and the  hypothesis imply that
$$0=k( u^2(e_0),e_0)=k(u(e_0),u(e_0))= x_1^2-x_2^2.$$
Then $x_1=\pm x_2.$  Suppose first that $x_1=x_2=a$ ($a\ne
0,$because $u(e_0)\notin \mathbb{R} e_0).$ Then $u(e_0)=\alpha e_0
+a(e_1+e_2).$ Consider a new basis of ${\cal{G}}$ consisting of
$e_{-1}, e_0, v_1=e_1+e_2, v_2=e_1-e_2.$ Then
$$
V:={\mathrm{ker}}\,{\mathrm{ad}}_{u(e_0)}={\mathrm{Span}}\,\{e_0,
v_1\}.
$$
We have that $u(e_{-1})=-x_2e_1-x_1e_2=-a v_1.$ Hence
$$
{\mathrm{im}}\,{\mathrm{ad}}_{u(e_0)}={\mathrm{Span}}\,\{e_0,
u(e_1)\}={\mathrm{Span}}\,\{e_0, v_1\}.
$$
This implies that $u$ leaves $V$ invariant. In particular $u(v_1)\in
V.$ Hence $u(v_1)=Ae_0+Bv_1,$ and
$$
A=k(u(v_1), e_{-1})=k(v_1, u(e_{-1}))=-ak(v_1,v_1)=0.
$$
Then $u(v_1)$ and $u(e_{-1})$ are linearly dependent, which is not
possible.

If $x_1=-x_2=a,$ then the same proof applies, with $v_2$ playing the
role of $v_1.$

Finally,  suppose that ${\mathrm{rank}}\,{\mathrm{ad}}_{u(e_0)}= 0,$
that is $u(e_0)\in\mathbb{R} e_0. $ Without loss of generality,
suppose that $u(e_0)=e_0.$

Some calculations show that
$$
k\left(\mathrm{R}_{e_{-1},e_1}e_{-1},e_2\right)=
\frac{b(-1+a+d)}{b^2+ad}
$$
where $b=k(u(e_2),e_1),$ $a=k(u(e_1),e_1),$ and $d=-k(u(e_2),e_2).$
The semi Riemannian metric being flat, either $b=0$ or $a+d=1.$ If
$a+d=1$ then
$$k\left(\mathrm{R}_{e_{-1},e_1}e_{-1},e_1\right)=1.$$
Since the former equality contradicts the hypothesis, $b=0.$ Using
this, we get
$$
0=k\left(\mathrm{R}_{e_{-1},e_1}e_{-1},e_1\right)=
\frac{1}{4ad}(a^2+2a(-1+d)-(-1+d)(1+3d)),
$$
and
$$
0=k\left(\mathrm{R}_{e_{-1},e_2}e_{-1},e_2\right)=
\frac{1}{4ad}(3a^2-2a(1+d)-(-1+d)^2).$$ Adding the two equations, we
get
$$
0=a^2-a-d(-1+d)=(a-d)(a+d-1).
$$
It is easy to check that neither $a=d$ nor $a+d=1$ satisfy
$$
3a^2-2a(1+d)-(-1+d)^2=0.
$$

In order to conclude the proof of Theorem \ref{theorem:dim4}, it is
easy to check that the linear map ${\cal{G}}\to
\mathrm{gl}\,({\cal{G}})$ given by
$$
x=x_{-1}e_{-1}+x_{0}e_{0}+x_{1}e_{1}+x_{2}e_{2}\mapsto
L_x=\dfrac{1}{2}\left(
\begin{array}{cccc}
0 & 0 & 0 & 0 \\
0 & 0 & x_2 &-x_1  \\
0 & 0 & 0 & 2x_{-1} \\
0 & 0 & 2x_{-1} & 0 \\
\end{array}
\right)
$$
is a left symmetric product on ${\cal{G}}$ compatible with the Lie
bracket. \qed

\begin{rmk}
In fact, for this latter quadratic group there is a left invariant
affine structure  which is holomorphic (\cite{kn:dm}).
\end{rmk}

As a consequence, we have that for undecomposable quadratic Lie
groups of dimension 4, no exponential map of a left invariant semi
Riemannian metric is a local isometry.

\vspace{.1in} The following results give flat left invariant semi
Riemannian metrics on nilpotent Lie groups. Let $f$ be an
endomorphism of the underlying vector space of a Lie algebra
${\cal{G}},$ such that $$
f([x,y])-[f(x),f(y)]\in{\cal{Z}}\,({\cal{G}})
$$
for all $x,y\in{\cal{G}}.$ Such an endomorphism is called a
{\bf{q-homomorphism}} of Lie algebras. An endomorphism $d$ of the
linear space ${\cal{G}}$ is called an $f$ derivation if
$$
d[x,y]=[dx,fy]+[fx,dy],
$$
for all $x, y\in {\cal{G}}.$ In particular a derivation is a
${\mathrm{Id}}$ derivation.

\begin{pr}
Let $({\cal{G}},k)$ a quadratic Lie algebra with an invertible $f$
derivation $d.$  Then the semi Riemannian metric defined by
$$\langle x,y\rangle =k(dx,dy)$$ is flat and the Levi-Civita product is given by
$$  x y =d^{-1}[fx,dy].$$
\end{pr}

\vspace{.1in}\noindent{\bf{Proof.}} The Levi-Civita product is given by the
Koszul formula:
$$
2\langle xy,z\rangle=\langle [x,y],z\rangle-\langle [y,z],x \rangle+\langle [z,x],y\rangle
$$
By the definition of the metric,
$$
2k(d(xy),dz)=k(d[x,y],dz)-k(d[y,z],dx)+ k(d[z,x],dy).
$$
Using the fact that $d$ is a $f$ derivation,
$$
\begin{array}{rcl}
2k(d(xy),dz)&=&k([dx,fy],dz)-k([dy,fz],dx)+ k([dz,fx],dy)\\
&&+k([fx,dy],dz)-k([fy,dz],dx)+ k([fz,dx],dy)\\
&=& 2k([fx,dy],dz).
\end{array}
$$
Hence
$$
d(xy)=[fx,dy]
$$
because $k$ is non degenerate.
A simple calculation shows that
$$
(xy)z=d^{-1}[fx,[fy,dz]].
$$
Therefore,
$$
\begin{array}{rcl}
(xy)z-(yx)z&=&d^{-1}\left([fx,[fy,dz]]-[fy,[fx,dz]]\right)\\
           &=&d^{-1}[[fx,fy],dz]
\end{array}.
$$
Finally,
$$
[x,y]z=d^{-1}\left([f[x,y],dz]\right)=d^{-1}\left([[fx,fy],dz]\right)
$$
because $d$ is a $f$ derivation and $f$ is a q-homomorphism of Lie algebras.

\begin{theorem}\label{theorem:fder3}
Every 3 step nilpotent Lie group has an invertible $f$ derivation
$d$ that induces a flat left
invariant connection on $G.$
\end{theorem}

\vspace{.2in} \noindent {\bf{Proof.}} Every  3 step nilpotent can be
decomposed as
$$
{\cal{G}}={\cal{G}}_0\oplus{\cal{G}}_1\oplus{\cal{G}}_2
$$
where ${\cal{G}}_0=[{\cal{G}},[{\cal{G}},{\cal{G}}]]\subset
{\cal{Z}}({\cal{G}}),$ ${\cal{G}}_1$ is a supplement of
${\cal{G}}_0$ in $[{\cal{G}},{\cal{G}}]$ and ${\cal{G}}_2$ is a
supplement of $[{\cal{G}},{\cal{G}}]$ in ${\cal{G}}.$ Define
$$
f(x)=a_ix\quad {\mbox{for}}\quad x\in{\cal{G}}_i,
$$
where $a_1=4/9,$ $a_2=2/3,$ and
$$
d(x)=\alpha_i x \quad {\mbox{for}}\quad x\in{\cal{G}}_i,
$$
with $\alpha_0=\alpha_1.$ The conditions on the parameters in order
that $d$ is an $f$ derivation are:
$$
\begin{array}{ccc}
\alpha_0\alpha_2&\ne&0\\
\alpha_0&=&(4/9)\alpha_2+\alpha_0a_0\\
\alpha_2&=&1/3.
\end{array}
$$
The $f$ derivation $d$ is invertible and the
product
$$
xy:=d^{-1}([fx,dy])
$$ is left symmetric hence defines a flat left invariant connection
on ${\cal{G}}.$

\begin{theorem}\label{theorem:fder}
Every quadratic  3 step nilpotent Lie group admits a flat left
invariant semi Riemannian metric induced by an invertible $f$
derivation on its Lie algebra.
\end{theorem}

Another general situation with flat left invariant complete semi
Riemannian metric is for quadratic Lie groups with a left invariant
symplectic form (\cite{kn:MSRI}).

\begin{pr}\label{pr:nilp}
Let $({\cal{G}},k)$ be a nilpotent, quadratic  Lie algebra and $u$ a
$k$ symmetric isomorphism  of the vector space underlying
${\cal{G}}.$ If $u$ preserves the descending central sequence
${\cal{G}}$ (and hence the ascending central sequence of
${\cal{G}}$) then the metric $k_u$ is complete. Moreover, the
solutions of the Euler equation (\ref{eq:geo}) are polynomial.
\end{pr}

The proposition follows from the following lemma.
\begin{lm}\label{lm:nilp}
If $({\cal{G}},k)$ is a nilpotent, quadratic  Lie algebra of degree
$m$ and $v\in{\mathrm{End}}({\cal{G}})$ preserves the descending
central sequence of ${\cal{G}}$, then the $m$th derivative of the
vector field given by
$$ \dot{x}=[x,v(x)]
$$
is zero.
\end{lm}

\vspace{.1in} \noindent {\bf{Proof.}} Let $t\mapsto\alpha(t)$ be a
curve in  ${\cal{G}}.$ Define $\beta(t):=[\alpha(t), v(\alpha(t))].$
Then
$$
\forall\, i\in\mathbb{N}\setminus\{0\},\;\;\;
\beta^{(i)}(t)=\sum_{j=0}^i\, C_j^i\,[\alpha^{(j)}(t),
v(\alpha^{(i-j)}(t))]
$$
and $\beta^{(i)}\in{\mathrm{C}}^{(i+1)}({\cal{G}}).$ Hence, if
$x:t\mapsto x(t)$ is a solution of (\ref{eq:geo}) and
$\beta(t):=[x(t), v(x(t))],$ then $\beta^{(i)}(t)=x^{(i+1)}(t),$
$\forall i\in\mathbb{N}.$ It follows that
$x^{(m)}\in{\mathrm{C}}^m({\cal{G}})=\{0\}.$ \qed

\vspace{.1in} \noindent {\bf{Proof of Proposition \ref{pr:nilp}.}}
Consider in Lemma \ref{lm:nilp} the vector field
$\dot{x}=[x,u^{-1}(x)].$ \qed

\begin{rmk} There are incomplete left invariant semi Riemannian metrics on
nilpotent quadratic  Lie groups.
\end{rmk}

The following is an example of this situation.

\subsubsection*{Example.}\label{ex:1} Let
${\cal{G}}={\mathrm{Span}}\{e_0,e_1,e_2,e_3,e_4\}$ with Lie bracket
$$
[e_4,e_1]=e_2;\;\; [e_4,e_2]=e_3;\;\; [e_1,e_2]=e_0,
$$
the non stated products are obtained either by antisymmetry or are
zero. For $x\in{\cal{G}},$ let $x=
x_0e_0+x_1e_1+x_2e_2+x_3e_3+x_4e_4.$ The Lie algebra ${\cal{G}}$ is
3 step nilpotent and for $k$ given by
$$k(x,x):=2(x_0\,x_4-x_1\,x_3)+x_2^2,$$
$({\cal{G}},k)$ is quadratic. Let $u\in{\mathrm{GL}}({\cal{G}})$
with matrix given in the basis ${\cal{B}}=\{e_0,e_1,e_2,e_3,e_4\},$
$$
{\mathrm{M}}_{\cal{B}}(u)=\left(
\begin{array}{ccccc}
0&1&0&0&0\\
1&0&0&0&0\\
0&0&1&0&0\\
0&0&0&0&-1\\
0&0&0&-1&0
\end{array}
\right).
$$
It is easy to check that $u$ is $k$ symmetric. Equation
(\ref{eq:geo}) is (in the same coordinate system)
\begin{eqnarray*}
\dot{x}_0&=&-x_2x_0+x_1x_2\\
\dot{x}_1&=&0\\
\dot{x}_2&=&x_4x_0+x_1x_3\\
\dot{x}_3&=&x_2x_4+x_2x_3\\
\dot{x}_4&=&0.
\end{eqnarray*}
The curve
$$
x_0(t)=\frac{-2}{(1+t)^2};\;\; x_1(t)=0;\;\;
x_2(t)=\frac{2}{1+t};\;\;x_3(t)=c(1+t)^2-1;\;\; x_4(t)=1
$$
is a non complete solution of Equation (\ref{eq:geo}), see
\cite{kn:B-M}. Hence the semi Riemannian metric defined by $u$ is
not flat.

Notice that the quadratic Lie algebra $({\cal{G}},k)$ given in the
example above is a 3 step nilpotent quadratic algebra. The flat
metric given by Theorem \ref{theorem:fder} is of signature $(2,3).$
This algebra is undecomposable.

\vspace{.2in} \noindent {\bf{Proof of the undecomposability.}}
Remark that ${\cal{Z}}({\cal{G}})$ is totally isotropic and that any
ideal of dimension 1 is central. Let ${\cal{I}}$ an ideal of
dimension 2. Since ${\cal{I}}\cap{\cal{Z}}({\cal{G}})\ne(0),$
${\cal{I}}={\mathrm{Span}}\{x,y\}$ where $x=x_0e_0+x_3e_3,$
$y=y_0e_0+y_1e_1+y_2e_2+y_3e_3+y_4e_4.$ The vectors
$$
\begin{array}{ccc}
x_0e_0+x_3e_3&& \\
y_2e_0-y_4e_2&=& [e_1,y]\\
-y_1e_0-y_4e_3&=&[e_2,y]\\
y_1e_2+y_2e_3&=&[e_4,y]\\
\end{array}
$$
are in ${\cal{I}}.$ Since we are assuming that the ideal is of
dimension 2, the matrix
$$
\left(\begin{array}{ccc}
x_0&0&x_3 \\
y_2&-y_4&0\\
-y_1&0&-y_4\\
0&y_1&y_2\\
\end{array}\right)
$$
has rank at most 2. This implies that
$$
\begin{array}{ccc}
y_4 (x_3 y_1 - x_0 y_4)&=&0\\
y_2 (x_3 y_1 - x_0 y_4)&=&0\\
y_1 (x_3 y_1 - x_0 y_4)&=&0
\end{array}
$$
If $x_3 y_1 - x_0 y_4\ne 0,$ $y_1=y_2=y_4=0$ and
$y\in{\cal{Z}}({\cal{G}}).$ If $x_3 y_1 - x_0 y_4=0,$ then
$$
\begin{array}{rcl}
k(x,x)&=&0\\
k(x,y)&=&k(x_0e_0+x_3e_3,y_0e_0+y_1e_1+y_2e_2+y_3e_3+y_4e_4)\\
&=&x_0y_4-x_3y_1=0.
\end{array}
$$
and the ideal ${\cal{I}}$ is degenerate. \qed

\section{Quadratic 2-step nilpotent Lie groups}

Let $({\cal{G}},k)$ be a quadratic 2-step nilpotent Lie algebra with
$0$ corank, that is such that
$[{\cal{G}},{\cal{G}}]={\cal{Z}}({\cal{G}}).$ Under this hypothesis
the Lie algebra $({\cal{G}},k)$ is isomorphic to a quadratic Lie
algebra $(V\oplus V^{*}, \theta, k)$ where
$V^{*}=[{\cal{G}},{\cal{G}}],$ $\theta\in\Lambda^3(V),$
${\mathrm{rank}}\,\theta={\mathrm{dim}}\,V,$ the bracket is given by
$$
[(x,\alpha),(y,\beta)]=(0,\theta(x,y,\cdot)),
$$
and
$$
k((x,\alpha),(y,\beta))=\alpha(y)+\beta(x).
$$
Let $\phi\in{\mathrm{Gl}}(V)$ and define $u:V\oplus V^{*}\to V\oplus
V^{*}$ by $u(x, \alpha):=(\phi (x),^t\phi(\alpha))=(\phi
(x),\alpha\circ\phi).$ It is easily verified that
$u\in{\mathrm{GL}}(V\oplus V^{*}).$ Denote by
$\langle\,,\rangle_\phi$ the bilinear form induced on $V\oplus
V^{*}$ via $k$ by $u$ (hence by $\delta$). Then
$\langle\,,\rangle_\phi$ is non degenerate and $u$ is
$\langle\,,\rangle_\phi$ symmetric. We have

\begin{theorem}\label{theorem:2n}
Let $(G,k)$ be a quadratic Lie group with Lie algebra
${\cal{G}}:=(V\oplus V^{*}, \theta, k),$ as above. Then for every
$\phi\in{\mathrm{Gl}}(V)$ the metric $\langle\,,\rangle_\phi$
defines a flat and geodesically complete semi Riemannian metric on
$G,$ and $(G,\langle\,,\rangle_\phi),$
$G,\langle\,,\rangle_{\phi'})$ are isometric if and only if
$\phi'=\psi^{-1}\phi\psi$ for some $\psi\in{\mathrm{Gl}}(V).$
Moreover, if ${\mathrm{V}}\ge 9$ there are infinitely many non
isometric such metrics.

\end{theorem}

\vspace{.1in}\noindent{\bf{Proof.}} The Levi-Civita product
associated to $\langle\,,\rangle_\phi$ is given by
$$
2ab:=2L_ab=[a,b]+u^{-1}([a,u(b)]+[b,u(a)]).
$$
In order to prove that $L_{[a,b]}=[L_a,L_b]$ notice that, since
$V^{*}=[{\cal{G}},{\cal{G}}]={\cal{Z}}({\cal{G}})$ and $u^{-1}$
invariant, then
$$
a(bc)=b(ac)=[a,b]c=0,$$ for all $a,b,c\in{\cal{G}}.$ Hence
$\langle\,,\rangle_\phi$ is flat, and geodesically complete because
${\cal{G}}$ is unimodular.

A straightforward calculation shows that $({\cal{G}},\theta,k,\phi)$
and $({\cal{G}},\theta,k,\phi')$ are isometric if and only if there
exists $\psi\in{\mathrm{Gl}}(V)$ such that
$\phi'=\psi^{-1}\phi\psi.$

Finally the Vinberg-Elashvili classification Theorem \cite{kn:V-E}
implies that there are infinitely many non degenerate and non
conjugate 3-linear forms on $V$ when ${\mathrm{dim}}\, V\ge 9.$
Consequently there are infinitely many non isometric flat left
invariant semi Riemannian metrics on $G.$

\qed

Theorem \ref{theorem:2n} can be used to construct flat compact semi
Riemannian nilmanifolds as is shown by the following example.

\subsubsection*{Example} Consider the Lie algebra $A_d$
with basis $\{ e_1,e_2,e_3,f_1,f_2\}$ and Lie bracket
$$
[e_1,e_2]=f_2,\quad [e_3,e_4]=f_2,\quad [e_1,e_3]=f_1,\quad
[e_2,e_4]=df_2
$$
where $d$ is a square free integer. It is clear that $A_d$ is a
2-step nilpotent Lie algebra of 0 corank. Moreover if $d\ne d'$ the
$\mathbb{Q}$ algebras $A_d$ and $A_{d'}$ are not isomorphic (see
\cite{kn:sch}). Hence the simply connected Lie group $G$ of Lie
algebra $\phantom{}^{*}tA_d:=A_d^{*}\rtimes_{coadj}A_d$ has
lattices. Consequently the manifold $M=\Gamma\backslash G,$ where
$\Gamma$ is a lattice, have many flat semi Riemannian metrics.

\vspace{.1in} For more details on 2-step nilpotent quadratic Lie
algebras, see \cite{kn:r}.

\section*{Acknowledgments}
The authors wish to thank Ph. Revoy  for the fruitful
discussions during the elaboration of this work.

\end{document}